\documentclass[review]{elsarticle}

\usepackage{amssymb}
\usepackage{amsmath}
\usepackage{amsthm}
\usepackage{cleveref}
\usepackage{esint}
\usepackage{color}
\crefname{theorem}{Theorem}{theorem}
\newtheorem{theorem}{Theorem}[section]

\numberwithin{equation}{section}

\newtheorem{lemma}[theorem]{Lemma}

\theoremstyle{definition}
\newtheorem{remark}[theorem]{Remark}

\begin{document}
\begin{frontmatter}

\title{Existence and boundary asymptotic behavior of large solutions of Hessian equations \tnoteref{t1}}

\tnotetext[t1]{This work was supported by NSFC grant 11671316.}
\author[rvt1]{Shanshan Ma \corref{cor1}}
\ead{mss5221@stu.xjtu.edu.cn}
\author[rvt1]{Dongsheng Li }
\ead{lidsh@mail.xjtu.edu.cn}

\cortext[cor1]{Corresponding author.}

\address[rvt1]{School of Mathematics and Statistics, Xi'an Jiaotong University, Xi'an 710049, China}

\begin{abstract}
In this paper, we establish the existence of large solutions of Hessian equations and obtain a new boundary asymptotic behavior of solutions.
\end{abstract}


\end{frontmatter}


\section{Introduction}
Let $f\in C^1(0,\infty)$ be positive and nondecreasing, and $b\in C^{1,1}(\overline{\Omega})$ be positive in $\Omega$, where $\Omega\subseteq \mathbb{R}^n(n\geq 2)$ is a bounded domain with boundary of class $C^2$. In this paper we investigate the following $k$-Hessian equation ($1\le k\le n$):
\begin{equation}\label{e1.1}
\left\{ {\begin{array}{*{20}{l}}
{S_k(D^2u) =\sigma_k(\lambda)= b(x)f(u)\ \ \mbox{in}\ \Omega,}\\
{u = \infty\ \ \ \ \ \ \ \ \ \ \ \ \ \ \ \ \ \ \ \ \ \ \ \ \ \ \ \ \ \ \ \mbox{on}\ \partial\Omega,}
\end{array}} \right.
\end{equation}
where $\lambda=(\lambda_1,\lambda_2,\cdots,\lambda_n)$ are the eigenvalues of the Hessian matrix $D^2u$ and
$$
\sigma_k(\lambda)=\sum\limits_{1 \le {i_1} <  \cdots< {i_k} \le n} {{\lambda _{{i_1}}} \cdots  {\lambda _{{i_k}}}}
$$
is the $k^{th}$ elementary symmetric function of $\lambda$. For completeness, we also set $\sigma_0(\lambda)=1$ and $\sigma_k(\lambda)=0$ for $k>n$. The boundary condition means $u(x)\rightarrow +\infty$ as $d(x)=\mbox{dist}(x,\partial\Omega)\rightarrow 0^+$.

To work in the realm of elliptic operators, we have to restrict the class of functions and domains. Following \cite{CNS3}, a function $u\in C^2(\Omega)$ is called a $k$-admissible function if  for any $x\in\Omega$, ${\lambda}(D^2u(x))$ belongs to the cone given by
$$
\Gamma_k=\{\lambda\in \mathbb{R}^n:\sigma_j(\lambda)>0, j=1,\cdots,k\}.
$$
From \cite{CNS3} and \cite{Gar}, $\Gamma_k$ is an open, convex, symmetric (under the interchange of any two $\lambda_j$) cone with vertex at the origin, and
$$\Gamma_1\supset\Gamma_2\supset\cdots\supset\Gamma_n=\{\lambda\in \mathbb{R}^n:\lambda_1,\cdots,\lambda_n>0\}.
$$
In \cite{CNS3}, it is also shown that
$$
\frac{\partial\sigma_k(\lambda)}{\partial\lambda_i}>0\ \ \mbox{in}\ \ \Gamma_k\ \ \forall\ i
$$
and
$$
\sigma_k^{{1}/{k}}(\lambda)\ \mbox{is a concave function in}\ \Gamma_k.
$$

Let
$$
S(\Gamma_k)=\{A: A\in \mathbb{S}^{n\times n}, \lambda(A)\in\Gamma_k\},
$$
where $\mathbb{S}^{n\times n}$ denotes the set of $n\times n$ real symmetric matrices. Recall that $S(\Gamma_k)$ is an open convex cone with vertex at the origin in matrix spaces.
The properties of $\sigma_k$ described above guarantee that
$$
\left(\frac{\partial S_k}{\partial A_{ij}}\right)_{n\times n}>0\ \ \ \forall\ A\in S(\Gamma_k)
$$
and
$$
S_k^{{1}/{k}}\ \mbox{is concave in}\ S(\Gamma_k).
$$

For an open bounded subset $\Omega\subset \mathbb{R}^n$ with boundary of class $C^2$ and for every $x\in\partial\Omega$, we denote by $\rho(x)=(\rho_1(x),\cdots,\rho_{n-1}(x))$ the principal curvatures of $\partial\Omega$ (relative to the interior normal). Recall that $\Omega$ is said to be $l$-convex ($1\le l\le n-1$) if $\partial\Omega$, regarded as a hypersurface in $\mathbb{R}^n$, is $l$-convex, that is, for every $x\in \partial\Omega$, $\sigma_j(\rho(x))\ge 0$ with $j=1,2,\cdots,l$. Respectively, $\Omega$ is called strictly $l$-convex if $\sigma_j(\rho(x))>0$ with $j=1,2,\cdots,l$.

Since we will consider viscosity solutions of (1.1), we first give the following definitions (See \cite{Urb}). A function $u\in C(\Omega)$ is said to be a viscosity subsolution (supersolution) of (1.1) if $x_0\in\Omega$, $A$ is an open neighborhood of $x_0$, $\psi\in C^2(A)$ is $k$-admissible and $u-\psi$ has a local maximum (minimum) at $x_0$, then
$$
S_k(D^2\psi(x_0))\ge b(x_0)f(\psi(x_0))~(\le b(x_0)f(\psi(x_0))).
$$
A function $u\in C(\Omega)$ is said to be a viscosity solution if it is both a viscosity subsolution and a viscosity supersolution.

From now on we shall always assume that $\Omega$ is bounded and strictly $(k-1)$-convex and consider viscosity solutions of (1.1).

~

In order to show the existence of the viscosity solution of (1.1), we need $f$ to satisfy some conditions. Specifically, we assume that $f$ satisfies:

$\mathbf{(f_1)}~f\in C^{1}(0,\infty),~f(s)>0,$ and is nondecreasing in $(0,\infty)$;

$\mathbf{(f_2)}$ The function
$$
\Phi(s)=\int\limits_s^{\infty} {\frac{{d\tau }}{{H(\tau )}}}
$$
is well defined for any $s>0$, where
$$
\ H(\tau) = {((k + 1)F(\tau))^{{1 \mathord{\left/
 {\vphantom {1 {(k + 1)}}} \right.
 \kern-\nulldelimiterspace} {(k + 1)}}}}\ \ \ \forall\ \tau>0
$$
and
$$
\ F(\tau ) = \int\limits_0^{\tau}  {f(s)ds}\ \ \ \forall\ \tau>0.
$$
Since we will investigate the boundary behavior of $u$ and $u=\infty$ on the boundary, we only need to concern the behavior of $F(\tau)$ and $f(\tau)$ as $\tau$ being sufficiently large.
For convenience, we define by $\varphi$ the inverse of $\Phi$, i.e., $\varphi$ satisfies
\begin{equation}\label{Defvarphi}
\int\limits_{\varphi (t)}^{\infty} {{\frac{d\tau}{H(\tau)}} }  = t\ \ \ \forall\ t>0.
\end{equation}

Our existence results are stated as follows:
\begin{theorem}
Let $\Omega\subseteq \mathbb{R}^n$ be a bounded and strictly $(k-1)$-convex domain with  $\partial\Omega\in C^{3,1}$. Suppose that $f$ satisfies $\mathbf{(f_1)}$ and $\mathbf{(f_2)}$, and that $b\in C^{1,1}(\overline{\Omega})$ is positive in $\Omega$. Then problem (1.1) admits a viscosity solution $u\in C(\Omega)$.
\end{theorem}

The existence of viscosity solutions of Hessian equations with boundary blowup has been considered in Salani \cite{Sal}, where Salani obtained the existence of solutions by using radial function to constructed barrier functions. In this paper, we present a different proof from that in \cite{Sal}. For more existence and nonexistence results, we refer to \cite{CSF, Hua, Jian, Sal} and the references therein.

~

To study the boundary behavior of the solution of (1.1), we need $f$ and $b$ to satisfy more conditions. Precisely, we assume that $f$ satisfies:

$\mathbf{(f_3)}$ There exists $C_f>0$ such that
\begin{equation}
\mathop {\lim }\limits_{s \to \infty} {H^\prime}(s)\int\limits_s^{\infty} {\frac{{d\tau }}{{H(\tau )}}}  =  {C_f},
\end{equation}
where $H(\tau)$ is defined in $\mathbf{(f_2)}$.

We also assume that $b$ satisfies:

$\mathbf{(b_1)}\ b\in C^{1,1}(\overline{\Omega})$ is positive in $\Omega$;

$\mathbf{(b_2)}\ $There exist a positive and nondecreasing function $m(t)\in C^1(0,\delta_0)$ (for some $\delta_0>0$),  and two positive constants $\overline{b}$ and $\underline{b}$ such that
$$
\underline{b} = \mathop {\lim \inf}\limits_{\mathop {{x} \in\Omega }\limits_{d(x) \to 0}}  \frac{{b(x)}}{{{m^{k + 1}}(d(x))}} \le \mathop {\lim \sup }\limits_{\mathop {{x} \in\Omega }\limits_{d(x) \to 0}} \frac{{b(x)}}{{{m^{k + 1}}(d(x))}} = \overline{b},
$$
where $d(x)=\mbox{dist}(x,\partial\Omega)$. Moreover, there exists $C_m\in[0,\infty)$ such that
$$
\mathop {\lim }\limits_{t \to 0^+} \left( {\frac{{M(t)}}{{m(t)}}} \right)^{\prime} = {C_m},
$$
where $M(t) = \int\limits_0^t {m(s)ds}<\infty$, $0<t<\delta_0$.

The boundary behavior of the solution of (1.1) may involve the curvatures of $\partial\Omega$.
We set
\begin{equation}\label{Ll}
L_0=\max\limits_{\overline{x}\in\partial\Omega}\sigma_{k-1}(\rho(\overline{x}))\ \ \mbox{and}\ \ l_0=\min\limits_{\overline{x}\in\partial\Omega}\sigma_{k-1}(\rho(\overline{x})),
\end{equation}
where $\rho(\overline{x})=(\rho_1(\overline{x}),\rho_2(\overline{x}),\cdots,\rho_{n-1}(\overline{x}))$ are the principal curvatures of $\partial\Omega$ at $\overline{x}$. Observe that $0<l_0\le L_0<+\infty$ since $\Omega$ is bounded and strictly $(k-1)$-convex. The boundary estimates of the solution of (1.1) are related to $L_0$ and $l_0$.

Now, we state our boundary behavior results as follows.
\begin{theorem}\label{tm1.2}
Let $\Omega\subseteq \mathbb{R}^n$ be a bounded and strictly $(k-1)$-convex domain with $\partial\Omega\in C^{3,1}$. Suppose that $f$ satisfies $\mathbf{(f_1)}$, $\mathbf{(f_2)}$ and $\mathbf{(f_3)}$, and that $b$ satisfies $\mathbf{(b_1)}$ and $\mathbf{(b_2)}$. If
\begin{equation}\label{CfCm}
C_f>1-C_m,
\end{equation}
where $C_f$ and $C_m$ are the constants defined in $\mathbf{(f_3)}$ and $\mathbf{(b_2)}$ respectively, then the viscosity solution $u$ of (1.1) satisfies
\begin{equation}\label{MBB}
 1\le\mathop {\liminf}\limits_{\mathop {{x} \in\Omega }\limits_{d(x) \to 0}}  \frac{{u(x)}}{{\varphi ({\overline{\xi} }M(d(x)))}}\ \ \mbox{and}\ \ \mathop {\limsup}\limits_{\mathop {{x} \in\Omega }\limits_{d(x) \to 0}} \frac{{u(x)}}{{\varphi ({\underline{\xi} }M(d(x)))}}\le 1,
\end{equation}
where $\varphi$ is defined by (\ref{Defvarphi}),
\begin{equation}
{\underline{\xi}} = {\left( {\frac{{{\underline{b}}}}{{{L_0}(1 - C_f^{ - 1}(1 - {C_m}))}}} \right)^{{1 \mathord{\left/
 {\vphantom {1 {(k + 1)}}} \right.
 \kern-\nulldelimiterspace} {(k + 1)}}}}\ \ and\ \ \ {\overline{\xi}} = {\left( {\frac{{{\overline{b}}}}{{{l_0}(1 - C_f^{ - 1}(1 - {C_m}))}}} \right)^{{1 \mathord{\left/
 {\vphantom {1 {(k + 1)}}} \right.
 \kern-\nulldelimiterspace} {(k + 1)}}}}.
\end{equation}
Here $L_0$ and $l_0$ are the constants in (\ref{Ll}).
\end{theorem}

\begin{remark}
Lemmas 2.1 and 2.2 of \cite{Zha} showed that if $b$ satisfies $\mathbf{(b_1)}$ and $\mathbf{(b_2)}$, and $f$ satisfies $\mathbf{(f_1)}$, $\mathbf{(f_2)}$ and $\mathbf{(f_3)}$, then $0\le C_m\le 1$ and $C_f\ge 1$. Hence (\ref{CfCm}) holds if $C_f>1$, or $C_f=1$ with $C_m>0$.
\end{remark}

\begin{remark} Theorem 1.2 contains the following interesting cases.

$(\mbox{i})$ $b\equiv 1$ and $f(s)=s^{\gamma}$, $\gamma>k$. In this case, choose $m(t)=1$, and we get
$$
M(t)=t~~\mbox{and}~~C_m=1.
$$
We also obtain $C_f=\frac{\gamma+1}{\gamma-k}$,
$$
\varphi(t)=\left(\frac{(k+1)^{k}(\gamma+1)} {(\gamma-k)^{k+1}}\right)^{1/(\gamma-k)}t^{-(k+1)/(\gamma-k)},
$$
$$
{\underline{\xi}} =\left(\frac{1}{L_0}\right)^{1/(k+1)}~~\mbox{and}~~~{\overline{\xi}} =\left(\frac{1}{l_0}\right)^{1/(k+1)}.
$$
Note that in this case (1.5) holds for any $\gamma>k$. Hence, by (1.6), the solution of (1.1) satisfies
$$
1 \le \mathop {\lim \inf}\limits_{\mathop {{x} \in\Omega }\limits_{d(x) \to 0}} \frac{{u(x)}}{{\left(\frac{l_0(k+1)^{k}(\gamma+1)} {(\gamma-k)^{k+1}}\right)^{1/(\gamma-k)}d(x)^{-(k+1)/(\gamma-k)}}}
$$
and
$$
\mathop {\lim \sup}\limits_{\mathop {{x} \in\Omega }\limits_{d(x) \to 0}} \frac{{u(x)}}{{\left(\frac{L_0(k+1)^{k}(\gamma+1)} {(\gamma-k)^{k+1}}\right)^{1/(\gamma-k)}d(x)^{-(k+1)/(\gamma-k)}}}\le1.
$$

$(\mbox{ii})$ $b= d(x)^{\alpha(k+1)}$, $\alpha>0$, near $\partial\Omega$ and $f(s)=s^{\gamma}$, $\gamma>k$. In this case, choose $m(t)=t^{\alpha}$, and we obtain
$$
M(t)=\frac{t^{\alpha+1}}{\alpha+1}~~\mbox{and}~~C_m=\frac{1}{\alpha+1}.
$$
We still have $C_f=\frac{\gamma+1}{\gamma-k}$,
$$
\varphi(t)=\left(\frac{(k+1)^{k}(\gamma+1)} {(\gamma-k)^{k+1}}\right)^{1/(\gamma-k)}t^{-(k+1)/(\gamma-k)},
$$
$$
{\underline{\xi}} = {\left( {\frac{{{(\alpha+1)(\gamma+1)}}}{{{L_0}({\gamma+1+\alpha k+\alpha})}}} \right)^{{1 \mathord{\left/
 {\vphantom {1 {(k + 1)}}} \right.
 \kern-\nulldelimiterspace} {(k + 1)}}}}~~\mbox{and}~~~{\overline{\xi}} = {\left( {\frac{{{(\alpha+1)(\gamma+1)}}}{{{l_0}({\gamma+1+\alpha k+\alpha})}}} \right)^{{1 \mathord{\left/
 {\vphantom {1 {(k + 1)}}} \right.
 \kern-\nulldelimiterspace} {(k + 1)}}}}.
$$
Note that in this case (1.5) holds for any $\gamma>k$. Therefore, by (1.6), the solution of (1.1) satisfies
$$
1 \le \mathop {\lim \inf}\limits_{\mathop {{x} \in\Omega }\limits_{d(x) \to 0}} \frac{{u(x)}}{{\left(\frac{l_0(\gamma+\alpha k+\alpha+1)(k+1)^k(\alpha+1)^k}{(\gamma-k)^{k+1}}\right)^{{1}/(\gamma-k)}d(x)^{-{(k+1)(\alpha+1)/(\gamma -k)}}}}
$$
and
$$
\mathop {\lim \sup}\limits_{\mathop {{x} \in\Omega }\limits_{d(x) \to 0}} \frac{{u(x)}}{{{\left(\frac{L_0(\gamma+\alpha k+\alpha+1)(k+1)^k(\alpha+1)^k}{(\gamma-k)^{k+1}}\right)^{{1}/(\gamma-k)}d(x)^{-{(k+1)(\alpha+1)/(\gamma -k)}}}}}\le1.
$$

\end{remark}

~

Problem (1.1) is the Laplace equation for $k=1$. The study of boundary blowup solutions of Laplace equation can be traced back to Bieberbach \cite{Bie}. The author considered $\Delta u=e^u$ in a smooth bounded domain in $\mathbb{R}^2$. Since then many papers have been dedicated to resolving existence, uniqueness and asymptotic behavior issues for solutions of blowup elliptic equations.  See \cite{BM92,BM95,CR02,CR03,G-M06,Kel,Oss} and their references.

For $k=n$, problem (1.1) is the Monge-Amp\`ere equation. There are many papers resolving existence, nonexistence, uniqueness and asymptotic behavior of boundary blowup solutions of Monge-Amp\`ere equation. We refer to \cite{CT,GJ,LM,LN,Moh,Zha} and the references therein.

For general $k$-Hessian equation with boundary blowup, there are also several authors studying the asymptotic behavior. In \cite{Sal}, Salani showed that: let $\Omega$ be a bounded and strictly convex domain, which implies that there exist two positive numbers $R_1\le R_2$ such that for any $y\in\partial\Omega$, there exist two balls $B_{1}^y$ and $B_{2}^y$, with rays $R_1$ and $R_2$ respectively, with $y\in\partial B_{1}^y\cap \partial B_{2}^y$ and $B_{1}^y\subset \Omega\subset B_{2}^y$. Suppose that $f$ satisfies $\mathbf{(f_1)}$ and $\mathbf{(f_2)}$, and $c_1\le b(x)\le c_2$ with $c_1$ and $c_2$ being two positive constants. Then the solution $u$ of (1.1) satisfies
$$
c_1^{1/(k+1)}p(R_1)\le\mathop {\lim \inf}\limits_{\mathop {{x} \in\Omega }\limits_{d(x) \to 0}}\frac{\Phi(u(x))}{d(x)}\le\mathop {\lim \sup}\limits_{\mathop {{x} \in\Omega }\limits_{d(x) \to 0}}\frac{\Phi(u(x))}{d(x)}\le c_2^{1/(k+1)}p(R_2),
$$
where $\Phi$ is the function in $\mathbf{(f_2)}$  and $p(R)=(C_{n-1}^{k-1})^{-1/(k+1)}R^{(k-1)/(k+1)}$.

Later, Huang \cite{Hua} generalized the results of \cite{Sal}. Let $\Omega$ be a smooth, strictly $(k-1)$-convex bounded domain. Suppose that~$f\in RV_q, q>k$, and $b$ satisfies $\mathbf{({b_1})}$ and $\mathbf{({b_2})}$. Define
$$
\mathcal{P}(\tau)=\sup\left\{\frac{f(y)}{y^k}: \theta\le y\le\tau\right\},~~~\mbox{for}~\tau\ge\theta,~ \theta~\mbox{sufficiently large},
$$
and
$$
\mathcal{P}^{\leftarrow}(s)=\inf\left\{\tau:\mathcal{P}(\tau)\ge s\right\}.
$$
Then the solution $u$ of (1.1) satisfies
$$
\xi^{-}\le \liminf\limits_{d(x)\rightarrow 0}{\frac{u}{\phi(d(x))}}~~\mbox{and}~~ \limsup\limits_{d(x)\rightarrow 0}{\frac{u}{\phi(d(x))}}\le \xi^{+},
$$
where,
$$
\phi(t)=\mathcal{P}^{\leftarrow}\left(\left(M(t)\right)^{-k-1}\right),~~~\mbox{for}~t>0~\mbox{small},
$$
and
$$
\frac{(\xi^{+})^{k-q}}{\underline{b}}\max\limits_{\partial\Omega}\sigma_{k-1}=\frac{(\xi^{-})^{k-q}}{\overline{b}}\min\limits_{\partial\Omega}\sigma_{k-1}=\frac{((q-k)/(k+1))^{k+1}}{1+C_m(q-k)/(k+1)}.
$$
Here $\underline{b}$, $\underline{b}$ and $C_m$ are given by $\mathbf{(b_2)}$.

In this paper, we investigate a new boundary behavior of solutions of (1.1).
In fact, we generalize the asymptotic results for Monge-Amp\`ere equation in \cite{Zha} to $k$-Hessian equation. Our results are also more accurate than \cite{Hua}.

Theorem 1.1 and 1.2 will be proved in Section 2 and 3 respectively.

\section{Proof of Theorem 1.1}
The following comparison principle is a basic tool for proofs of both Theorem 1.1 and Theorem 1.2.
\begin{lemma}[The comparison principle] Let $\Omega$ be a bounded domain in $\mathbb{R}^n$. Suppose that $g(x,\eta)$ is positive and continuously differentiable, and is nondecreasing only with respect to $\eta$. If $u,v\in C({\Omega})$ are respectively viscosity subsolution and supersolution of
$$
S_k(D^2u)= g(x,u)
$$
and $u\leq v$ on $\partial\Omega$,
then we have
$$
u\leq v\ \ \ \mbox{in}\ \Omega.
$$
\end{lemma}
\noindent\textbf{Proof.} We refer to Proposition 2.3 of \cite{Urb} for the detailed proof.

\qed

~

From \cite{Ivo} and \cite{IF}, we see, for any symmetric matrix $S$,
$$
S_k(S + \xi\times\xi) = S_k(S) + \frac{\partial S_k(S)}{\partial S_{ij}}\xi_i\xi_j, ~1\le k \le n,~\xi\in\mathbb{R}^n.
$$
Similarly,
$$
S_k(S - \xi\times\xi) = S_k(S) - \frac{\partial S_k(S)}{\partial S_{ij}}\xi_i\xi_j, ~1\le k \le n,~\xi\in\mathbb{R}^n.
$$
Then the following conclusion about composite functions can be obtained:
\begin{lemma} Let $\Omega$ be a bounded domain in $\mathbb{R}^n$. Suppose $h\in C^2(\mathbb{R})$ and $g\in C^2(\Omega)$. Then
$$
S_k(D^2h(g))=(h'(g))^{k-1}h''(g)S_k^{ij}(D^2g)g_ig_j+(h'(g))^kS_k(D^2g)~~~in~\Omega,
$$
where $S_k^{ij}(D^2g)=\frac{\partial S_k(D^2g)}{\partial g_{ij}}$.
\end{lemma}

Using the above lemma, we prove the following conclusion:

\begin{lemma} Let~$\Omega\subseteq \mathbb{R}^n$~be a bounded and strictly ~$(k-1)$-convex domain with $\partial\Omega\in C^{3,1}$.
Suppose that $f$ satisfies $\mathbf{(f_1)}$ and $\mathbf{({f_2})}$, and that $b\in C^{1,1}(\overline{\Omega})$~is positive. Then there exists a $\overline{h}\in C^2(\Omega)$, $\overline{h}(x)\rightarrow \infty$ as $\mbox{dist}(x,\partial\Omega)\rightarrow 0$, such that for any $k$-admissible function $u\in C^2(\Omega)\cap C(\overline{\Omega})$ satisfying
$$
S_k(D^2u)=b(x)f(u)~~~in~\Omega,
$$
we have
$$
u\le \overline{h}~~~in~\Omega.
$$
\end{lemma}

\noindent\textbf{Proof.} We assume that $w(w<0)$ is the admissible solution of
\begin{equation}\label{rightb}
\left\{ {\begin{aligned}
&S_k(D^2 w) = b(x) &&~~~\mbox{in}~\Omega;\\
&w = 0 &&~~~\mbox{on}~\partial\Omega,
\end{aligned}} \right.
\end{equation}
where $b\in C^{1,1}(\overline{\Omega})$ is positive. Indeed, from Theorem 1.1 of \cite{Tru95}, (\ref{rightb}) is uniquely solvable for admissible $w\in C^{3,\beta}(\overline{\Omega})$ for any $0<\beta<1$.

Define
$$
\overline{h}=\varphi(-\varepsilon w)~~~\mbox{in}~\Omega,
$$
where $\varphi$ is defined by (\ref{Defvarphi}) and $\varepsilon>0$ is a constant to be chosen later, and we take the second derivative of $\overline{h}$,
\begin{align*}
\overline{h}_{ij}&=(\varphi(-\varepsilon w))_{ij}\\
&=-\varepsilon\varphi{'}(-\varepsilon w)w_{ij}+\varepsilon^2\varphi{''}(-\varepsilon w)w_iw_j.
\end{align*}
By lemma 2.2, for $1\le l\le k$,
$$
\begin{array}{l}
S_l(D^2\overline{h})\\
=\varepsilon^l(-\varphi{'}(-\varepsilon w))^lS_l(D^2w)+\varepsilon^{l+1}\varphi{''}(-\varepsilon w)(-\varphi{'}(-\varepsilon w))^{l-1}S_{l}^{ij}(D^2w)w_iw_j\\
> 0~~~\mbox{in}~\Omega
\end{array}
$$
since $-\varphi'>0,~\varphi''\ge0$ and $w$ is $k$-admissible. That is $\overline{h}$ is $k$-admissible.
Specially,
\begin{equation}\label{SkOh}
\begin{array}{l}
S_k(D^2\overline{h})\\
=\varepsilon^k(-\varphi{'}(-\varepsilon w))^kS_k(D^2w)+\varepsilon^{k+1}\varphi{''}(-\varepsilon w)(-\varphi{'}(-\varepsilon w))^{k-1}S_{k}^{ij}(D^2w)w_iw_j\\
=\varepsilon^k\varphi{''}(-\varepsilon w)(-\varphi{'}(-\varepsilon w))^{k-1}S_k(D^2w)\left[\frac{-\varphi{'}(-\varepsilon w)}{\varphi{''}(-\varepsilon w)}+\varepsilon\frac{S_{k}^{ij}(D^2w)w_iw_j}{S_k(D^2w)}\right]\\
=\varepsilon^kb(x)f(\overline{h})\left[\frac{-\varphi{'}(-\varepsilon w)}{\varphi{''}(-\varepsilon w)}+\varepsilon\frac{S_{k}^{ij}(D^2w)w_iw_j}{S_k(D^2w)}\right].
\end{array}
\end{equation}
Let
$$
M_{\varepsilon}(x)=\varepsilon^k\left[\frac{-\varphi{'}(-\varepsilon w)}{\varphi{''}(-\varepsilon w)}+\varepsilon\frac{S_{k}^{ij}(D^2w)w_iw_j}{S_k(D^2w)}\right].
$$
We claim: for any $x\in\Omega$, $M_{\varepsilon}(x)$ is sufficiently small as $\varepsilon$ being sufficiently small.

In fact, by the choice of $\varphi$, we have
$$
\frac{-\varphi{'}(-\varepsilon w)}{\varphi{''}(-\varepsilon w)}=\frac{{(k+1)F(\varphi(-\varepsilon w) )}^{{{ k} \mathord{\left/
 {\vphantom {{ k} {k+1 }}} \right.
 \kern-\nulldelimiterspace} {k+1 }}}}{f(\varphi(-\varepsilon w))}.
$$
From \cite{GP} and \cite{Moh}, we know that if $f$ satisfies $\mathbf{({f_2})}$, then
\begin{equation}\label{Ff}
\mathop {\lim }\limits_{t \to \infty}\frac{{F(t )}^{{{ k} \mathord{\left/
 {\vphantom {{ k} {k+1 }}} \right.
 \kern-\nulldelimiterspace} {k+1 }}}}{f(t)}=0.
\end{equation}
Note that $\varphi(0)=\mathop {\lim }\limits_{t \to {0^+ }}\varphi(t)=+\infty$, we see
\begin{equation}\label{lim2.1.1}
\mathop {\lim }\limits_{\varepsilon \to 0^+}\frac{-\varphi{'}(-\varepsilon w)}{\varphi{''}(-\varepsilon w)}=0.
\end{equation}
Furthermore, since $w\in C^{3,\beta}(\overline{\Omega})$~is~$k$-admissible and the matrix $\{{S_{k}^{ij}(D^2w)}\}>0$, we have,
$$
\frac{S_{k}^{ij}(D^2w)w_iw_j}{S_k(D^2w)}>0~\mbox{is bounded in}~\Omega.
$$
This combining with (\ref{lim2.1.1}) imply that our claim holds.

From (\ref{SkOh}), we have, for sufficiently small $\varepsilon$,
$$
S_k(D^2{\overline{h}})\le b(x)f(\overline{h})~~~ \mbox{in}~\Omega.
$$
Note that, by the definition of $\overline{h}$, $\overline{h}(x)\rightarrow \infty$ as $\mbox{dist}(x,\partial\Omega)\rightarrow 0$. Therefore, for any $u\in C^2(\Omega)\cap C(\overline{\Omega})$ satsifying
$$
S_k(D^2u)=b(x)f(u)~~~\mbox{in}~\Omega,
$$
we have $u\le \overline{h} ~\mbox{on}~\partial\Omega$ and then by Lemma 2.1,
$$
u\le \overline{h}~~~\mbox{in}~\Omega.
$$
\qed

~

\noindent\textbf{Proof of Theorem 1.1.}
Since (\ref{Ff}) holds, we see, by $\mathbf{({f_2})}$, for sufficiently large $s>0$,
$$
\Psi(s)=\int\limits_s^{\infty} {\frac{1}{{f(\tau )}^{{{ 1} \mathord{\left/
 {\vphantom {{ 1} {k }}} \right.
 \kern-\nulldelimiterspace} {k }}}}d\tau }
$$
is well defined. Let~$\psi$~be the inverse of $\Psi$, i.e., $\psi$ satisfies
$$
s=\int\limits_{\psi(s)}^{\infty} {\frac{1}{{f(\tau )}^{{{ 1} \mathord{\left/
 {\vphantom {{ 1} {k }}} \right.
 \kern-\nulldelimiterspace} {k }}}}d\tau }.
$$
Then we have
$$
\psi(0)=\mathop {\lim }\limits_{t \to 0^+}=\infty,~\psi{'}(s)=-({f(\psi(s) )})^{{{ 1} \mathord{\left/
 {\vphantom {{ 1} {k }}} \right.
 \kern-\nulldelimiterspace} {k }}},~\psi{''}(s)=\frac{1}{k}{(f(\psi(s) ))}^{{{ (2-k)} \mathord{\left/
 {\vphantom {{ (2-k)} {k }}} \right.
 \kern-\nulldelimiterspace} {k }}}f{'}(\psi(s)).
$$

Assume that $w$ is the admissible solution of (\ref{rightb}) with $b\in C^{1,1}(\overline{\Omega})$ being positive in $\Omega$ (See \cite{Tru95}). We define
$$
\underline{h}(x)=\psi(-w(x)),\ \ \ {x\in}\ \Omega,
$$
and for $j=1,2,\cdots$,
$$
\Omega_j=\{x\in\Omega: \underline{h}(x)<{j}\}.
$$
Since~$w$~is~$k$-admissible, we see that $\Omega_j$ is strictly $(k-1)$-convex (see \cite{Tru97}).

Consider
\begin{equation}\label{EquinOmegaj}
\left\{ {\begin{aligned}
&S_k(D^2u) = b(x)f(u),&&~~~{x\in}~\Omega_j;\\
&u = {j},&&~~~{x\in}~\partial\Omega_j.
\end{aligned}} \right.
\end{equation}
We show that (\ref{EquinOmegaj}) has a $k$-admissible solution $u_j$. By Theorem 4.1 of \cite{Li}, we only need to prove that (\ref{EquinOmegaj}) has a $k$-admissible subsolution, and we will show that $\underline{h}$ is actually a $k$-admissible subsolution of (\ref{EquinOmegaj}).

By direct computation,
\begin{align*}
\underline{h}_{ij}&=(\psi(-w))_{ij}\\
&=-\psi{'}(-w)w_{ij}+\psi{''}(-w)w_iw_j.
\end{align*}
Since $f$~is positive and nondecreasing, $\psi^{\prime}<0$ and $\psi^{\prime\prime}\ge0$. Then according to that the matrix $\{w_iw_j\}$ is nonnegative, we have
$$
D^2\underline{h}\ge
-\psi{'}(-w) D^2w.
$$
This implies that for any $1\le j\le k$,
$$
S_j(D^2\underline{h})\ge S_j(-\psi{'}(-w) D^2w)>0,\ \ \ {x\in}\ \Omega,
$$
i.e., $\underline{h}$~is~$k$-admissible, and
$$
S_k(D^2\underline{h})\ge (-\psi{'}(-w))^kS_k(D^2w)=(-\psi{'}(-w))^kb(x)=b(x)f(\underline{h}),\ \ \ {x\in}\ \Omega
$$
since $w$ is $k$-admissible and satisfies (\ref{rightb}).
By the construction of $\Omega_j$,
$$
\underline{h}={j},\ \ \ {x\in}\ \partial\Omega_j.
$$
That is, we have shown that $\underline{h}$~is a $k$-admissible subsolution of (\ref{EquinOmegaj}). Therefore, (\ref{EquinOmegaj}) has a $k$-admissible solution $u_j$.

Since, for any $j=1,2,\cdots$,
$$
\underline{h}=u_j,\ \ \ {x\in}\ \partial\Omega_j,
$$
we have, by Lemma 2.1,
\begin{equation}\label{under h le uj}
\underline{h}\le u_j,~~~{x\in} ~\Omega_j.
\end{equation}
Furthermore, since
$$
u_j=\underline{h}\le u_{j+1},\ \ \ {x\in}\ \partial\Omega_j,
$$
we see,
\begin{equation}\label{mon inc}
u_j\le u_{j+1},~~~{x\in}~\Omega_j.
\end{equation}

For any $\Omega^{'}\subset\subset\Omega$ with $\Omega^{'}$ being strictly $(k-1)$-convex and $\partial\Omega^{'}\in C^{3,1}$, by Lemma 2.3, there exists a $\overline{h}\in C^2(\Omega')$, $\overline{h}(x)\rightarrow \infty$ as $\mbox{dist}(x,\partial\Omega')\rightarrow 0$, such that for sufficiently large $j$, $\Omega^{'}\subset\Omega_j$, and
$$
u_j\le \overline{h},~~~x\in~\Omega'.
$$
This combining with (\ref{mon inc}) imply that for any $x\in\Omega'$, the limit function
$$
u(x)=\mathop {\lim }\limits_{j \to {\infty }}u_j(x)
$$
exists. Then by the diagonal rule, for any $x\in \Omega$, $u(x)$ exists.
Moreover, $u(x)\rightarrow \infty$ as $d(x,\partial\Omega)\rightarrow 0$ and is a viscosity solution of (1.1).

\qed

\section{Proof of Theorem 1.2}
To study the boundary behavior of solutions of (1.1), we need the asymptotic estimate of functions in $\mathbf{(f_2)}$ and $\mathbf{(b_2)}$ as $t\rightarrow 0$. The following two lemmas describe those asymptotic behaviors.
\begin{lemma}
Let $m$ and $M$ be the functions given by $\mathbf{(b_2)}$. Then
$$
M(0)=\mathop {\lim }\limits_{t \to {0^+ }}M(t)=0,
$$
$$
\mathop {\lim }\limits_{t \to {0^+ }} \frac{{M\left( t \right)}}{{m\left( t \right)}} = 0,
$$
and
$$
\mathop {\lim }\limits_{t \to {0^+ }} \frac{{M\left( t \right)}m^{\prime}(t)}{{m^2\left( t \right)}}= 1-C_m.
$$
\end{lemma}

\begin{lemma}
Assume that $f$ satisfies $\mathbf{(f_1)},\ \mathbf{(f_2)}\ \mbox{and}\ \mathbf{(f_3)}$, and $\varphi$ satisfies (\ref{Defvarphi}). Then we have\\
$(i_1)\ \varphi(t)>0,\  \varphi(0)=\mathop {\lim }\limits_{t \to {0^+ }}\varphi(t)=+\infty,\ \varphi^{\prime}(t)=-\left((k+1)F(\varphi(t))\right)^{1 \mathord{\left/
 {\vphantom {1 {\left( {k + 1} \right)}}} \right.
 \kern-\nulldelimiterspace} {\left( {k + 1} \right)}},\\
 \indent{\ and}\ \varphi^{\prime\prime}(t)=\left((k+1)F(\varphi(t))\right)^{(1-k) \mathord{\left/
 {\vphantom {1 {\left( {k + 1} \right)}}} \right.
 \kern-\nulldelimiterspace} {\left( {k + 1} \right)}}f(\varphi(t));$\\
$(i_2)\mathop {\lim }\limits_{t \to {0^+ }}\frac{-\varphi^{\prime}(t)}{t\varphi^{\prime\prime}(t)}=\mathop {\lim }\limits_{t \to {0^+ }}\frac{{{{\left( {\left( {k + 1} \right)F\left( {\varphi \left( t \right)} \right)} \right)}^{{k \mathord{\left/
 {\vphantom {k {(k + 1)}}} \right.
 \kern-\nulldelimiterspace} {(k + 1)}}}}}}{{tf\left( {\varphi \left( t \right)} \right)}}=\frac{1}{C_f}.$
\end{lemma}

For detailed proofs of Lemmas 3.1 and 3.2, we refer to Lemmas 2.1 and 2.3 of \cite{Zha} respectively. More characterization of functions in $\mathbf{(f_2)}$ and $\mathbf{(b_2)}$ are also provided there.

~

We also need to recall some results of the distance function. Let $d(x)=\mbox{dist}(x,\partial\Omega)=\mathop {\inf }\limits_{{\rm{y}} \in \partial \Omega } \left| {x - y} \right|$. For any $\delta>0$, we define
\begin{equation}\label{DD}
\Omega_{\delta}=\{x\in\Omega:0<d(x)<\delta\}.
\end{equation}
If $\Omega$ is bounded and $\partial\Omega\in C^2$, by Lemma 14.16 of \cite{GT}, there exists $\delta_1>0$ such that
$$
d\in C^2(\Omega_{\delta_1}).
$$
Let $\overline{x}\in\partial\Omega$, satisfying $\mbox{dist}(x,\partial\Omega)=\left| {x - \overline{x}} \right|$, be the projection of the point $x\in\Omega_{\delta_1}$ to $\partial\Omega$, and $\rho_i(\overline{x})(i=1,\cdots,n-1)$ be the principal curvatures of $\partial\Omega$ at $\overline{x}$. Then, in terms of a principal coordinate system at $\overline{x}$, we have, by Lemma 14.17 of \cite{GT},
\begin{equation}\label{dD}
\left\{ {\begin{array}{l}
{Dd(x) = (0,0, \cdots ,1)},\\
{{D^2}d(x) = {\rm{diag}}\left[ {\frac{{ - {\rho _1}(\bar x)}}{{1 - d(x){\rho _1}(\bar x)}}, \cdots ,\frac{{ - {\rho _{n - 1}}(\bar x)}}{{1 - d(x){\rho _{n - 1}}(\bar x)}},0} \right].}
\end{array}} \right.
\end{equation}

~

\noindent\textbf{Proof of Theorem 1.1.} Under the assumptions in Theorem 1.2, we have the following conclusions: For any $\varepsilon>0$, we choose $\delta_{\varepsilon}>0$ small enough such that

$\mathbf{(a_1)}$ $m(t)$ satisfies $\mathbf{(b_2)}$ for $0<t<\delta_{\varepsilon}$;

$\mathbf{(a_2)}$ $d(x)\in C^2(\Omega_{2\delta_{\varepsilon}})$, where $\Omega_{2\delta_{\varepsilon}}$ is defined by (\ref{DD});

$\mathbf{(a_3)}$ $(\underline{b}-\varepsilon)m^{k+1}(d(x))\le b(x)\le(\overline{b}+\varepsilon)m^{k+1}(d(x))$ in $\Omega_{2\delta_{\varepsilon}}$;

$\mathbf{(a_4)}$ For any $1\le j\le k-1$, $\sigma_j\bigg(\frac{\rho_1(\overline{x})}{1-d(x)\rho_1(\overline{x})},\cdots,\frac{\rho_{n-1}(\overline{x})}{1-d(x)\rho_{n-1}(\overline{x})}\bigg)>0$ in $\Omega_{2\delta_{\varepsilon}}$.  Recall that $\rho_i(\overline{x})$ ($i=1,2,\cdots,n-1$) denote the principal curvatures of $\partial\Omega$ at $\overline{x}$, where $\overline{x}\in\partial\Omega$ satisfies ${d}(x)=\left| {x - \overline{x}} \right|$;

$\mathbf{(a_5)}$ $(1-\varepsilon)l_0\le \sigma_{k-1}\bigg(\frac{\rho_1(\overline{x})}{1-d(x)\rho_1(\overline{x})},\cdots,\frac{\rho_{n-1}(\overline{x})}{1-d(x)\rho_{n-1}(\overline{x})}\bigg)\le (1+\varepsilon)L_0$ in $\Omega_{2\delta_{\varepsilon}}$;

$\mathbf{(a_6)}$ $\sigma_k\bigg(\frac{\rho_1(\overline{x})}{1-d(x)\rho_1(\overline{x})},\cdots,\frac{\rho_{n-1}(\overline{x})}{1-d(x)\rho_{n-1}(\overline{x})}\bigg)$ is bounded in $\Omega_{2\delta_{\varepsilon}}$.

Fix $0<\varepsilon<\underline{b}/2$ and we choose
\begin{equation}\label{LConst}
{\underline{\xi} _{\varepsilon}} = {\left( {\frac{{{\underline{b}-2\varepsilon}}}{{{(1+\varepsilon)L_0}(1 - C_f^{ - 1}(1 - {C_m}))}}} \right)^{{1 \mathord{\left/
 {\vphantom {1 {(k + 1)}}} \right.
 \kern-\nulldelimiterspace} {(k + 1)}}}},
\end{equation}
and
\begin{equation}\label{UConst}
{\overline{\xi} _{\varepsilon}} = {\left( {\frac{{{\overline{b}+2\varepsilon}}}{{{(1-\varepsilon)l_0}(1 - C_f^{ - 1}(1 - {C_m}))}}} \right)^{{1 \mathord{\left/
 {\vphantom {1 {(k + 1)}}} \right.
 \kern-\nulldelimiterspace} {(k + 1)}}}},
\end{equation}
where $\underline{b}$, $\overline{b}$, $C_m$, $L_0$, $l_0$ and $C_f$ are given by $\mathbf{(b_2)}$, (\ref{Ll}) and $\mathbf{(f_3)}$ respectively. Let $\delta_{\varepsilon}$ be a small enough constant such that the above $\mathbf{(a_1-a_6)}$ hold and choose $0<\sigma<\delta_{\varepsilon}$. We define
\begin{equation}\label{d1xd2x}
d_1(x)=d(x)-\sigma, ~~d_2(x)=d(x)+\sigma
\end{equation}
and
\begin{equation}\label{supersub}
\left\{ {\begin{array}{l}
{\overline{u}}_{\varepsilon}(x)=\varphi(\underline{\xi} _{\varepsilon}M(d_1(x)))\ \ \mbox{in}\ \Omega_{2\delta_{\varepsilon}}/\overline{\Omega}_{\sigma},\\
{\underline{u}}_{\varepsilon}(x)=\varphi(\overline{\xi} _{\varepsilon}M(d_2(x)))\ \ \mbox{in}\ \Omega_{2\delta_{\varepsilon}-\sigma}.
\end{array}} \right.
\end{equation}

We divide the proof of Theorem 1.2 into three steps.

\textbf{Step 1.} We prove that ${\overline{u}}_{\varepsilon}$ is $k$-admissible and
\begin{equation}\label{supersolution}
S_k({D^2}{{\overline u} _\varepsilon }\left( x \right)) \leq b(x)f\left( {  {{{\overline u} }_\varepsilon }\left( x \right)}\right)\ \ \mbox{in}\ \Omega_{2\delta_\varepsilon}/\overline{\Omega}_{\sigma}
\end{equation}
as $\delta_{\varepsilon}$ being sufficiently small.

First, we show that ${\overline{u}}_{\varepsilon}$ is a $k$-admissible function in $\Omega_{2\delta_{\varepsilon}}/\overline{\Omega}_{\sigma}$. That is, for $1\le j\le k$,
\begin{equation}\label{superadmi}
S_j(D^2{\overline{u}}_{\varepsilon})>0~~\mbox{in}~~\Omega_{2\delta_{\varepsilon}}/\overline{\Omega}_{\sigma}.
\end{equation}

In view of (\ref{supersub}), we see, obviously,
$$
\overline{u}_{\varepsilon}(x)>0\ \ \mbox{in}\ \Omega_{2\delta_{\varepsilon}}/\overline{\Omega}_{\sigma}\ \ \mbox{and}\ \ \overline{u}_{\varepsilon}(x)=\infty\ \ \mbox{on}\ \partial{\Omega}_{\sigma}.
$$
By straightforward  computations,
$$
\begin{array}{l}
{\left( {{{\overline u }_\varepsilon }\left( x \right)} \right)_{\alpha\beta}}=(\varphi(\underline{\xi} _{\varepsilon}M(d_1(x))))_{\alpha\beta}\\
= {\underline{\xi} _{\varepsilon}}\left[ { {\underline{\xi} _{\varepsilon}}\varphi ''\left( {{\underline{\xi} _{\varepsilon}}M\left( {d_1\left( x \right)} \right)} \right){m^2}\left( {d_1\left( x \right)} \right)}
+\varphi '\left( {{\underline{\xi} _{\varepsilon}}M\left( {d_1\left( x \right)} \right)} \right)m'\left( {d_1\left( x \right)} \right)\right]{d_{\alpha}}{d_{\beta}}~~\\
~~ +{\underline{\xi} _{\varepsilon}}\ \varphi '\left( {{\underline{\xi} _{\varepsilon}}M\left( {d_1\left( x \right)} \right)} \right)m\left( {d_1\left( x \right)} \right){d_{\alpha\beta}}~
.
\end{array}
$$
Using (\ref{dD}) and Lemma 3.2$(i_1)$, we derive that for $1\le j\le k$,
\begin{equation}\label{superderi}
\begin{array}{l}
S_j(D^2{\overline{u}}_{\varepsilon})\\
={\underline{\xi}^{j+1}_{\varepsilon}}m^{j+1}(d_1(x))f(\varphi(\underline{\xi} _{\varepsilon}M(d_1(x)))){{{\left( {\left( {k + 1} \right)F\left( {\varphi \left( {{\underline{\xi} _{\varepsilon}}M\left( {{d_1}\left( x \right)} \right)} \right)} \right)} \right)}^{{(j-k) \mathord{\left/
 {\vphantom {(j-k) {\left( {k + 1} \right)}}} \right.
 \kern-\nulldelimiterspace} {\left( {k + 1} \right)}}}}}\\
~~\times\bigg[\left(1-{ \frac{{M\left( {{d_1}\left( x \right)} \right)m'\left( {{d_1}\left( x \right)} \right)}}{{{m^2}\left( {{d_1}\left( x \right)} \right)}}\frac{{{{\left( {\left( {k + 1} \right)F\left( {\varphi \left( {{\underline{\xi} _{\varepsilon}}M\left( {{d_1}\left( x \right)} \right)} \right)} \right)} \right)}^{{k \mathord{\left/
 {\vphantom {k {\left( {k + 1} \right)}}} \right.
 \kern-\nulldelimiterspace} {\left( {k + 1} \right)}}}}}}{{{\underline{\xi} _{\varepsilon}}M\left( {{d_1}\left( x \right)} \right)f\left( {\varphi \left( {{\underline{\xi} _{\varepsilon}}M\left( {{d_1}\left( x \right)} \right)} \right)} \right)}}}\right)\\
 ~~~~~~\times\sigma_{j-1}\left( \frac{{{\rho _1}\left( {\overline x } \right)}}{{1 - d\left( x \right){\rho _1}\left( {\overline x } \right)}},\cdots,\frac{{{\rho _{n-1}}\left( {\overline x } \right)}}{{1 - d\left( x \right){\rho _{n-1}}\left( {\overline x } \right)}}\right)\\
 ~~~~~~+\frac{{M\left( {{d_1}\left( x \right)} \right)}}{{{m}\left( {{d_1}\left( x \right)} \right)}}\frac{{{{\left( {\left( {k + 1} \right)F\left( {\varphi \left( {{\underline{\xi} _{\varepsilon}}M\left( {{d_1}\left( x \right)} \right)} \right)} \right)} \right)}^{{k \mathord{\left/
 {\vphantom {k {\left( {k + 1} \right)}}} \right.
 \kern-\nulldelimiterspace} {\left( {k + 1} \right)}}}}}}{{{\underline{\xi} _{\varepsilon}}M\left( {{d_1}\left( x \right)} \right)f\left( {\varphi \left( {{\underline{\xi} _{\varepsilon}}M\left( {{d_1}\left( x \right)} \right)} \right)} \right)}}\sigma_{j}\left( \frac{{{\rho _1}\left( {\overline x } \right)}}{{1 - d\left( x \right){\rho _1}\left( {\overline x } \right)}},\cdots,\frac{{{\rho _{n-1}}\left( {\overline x } \right)}}{{1 - d\left( x \right){\rho _{n-1}}\left( {\overline x } \right)}}\right)\bigg].
\end{array}
\end{equation}
This suggests us that to show (\ref{superadmi}), we only need to prove that for $1\le j\le k$,
\begin{equation}\label{impud}
\begin{array}{l}
\left(1-{ \frac{{M\left( {{d_1}\left( x \right)} \right)m'\left( {{d_1}\left( x \right)} \right)}}{{{m^2}\left( {{d_1}\left( x \right)} \right)}}\frac{{{{\left( {\left( {k + 1} \right)F\left( {\varphi \left( {{\underline{\xi} _{\varepsilon}}M\left( {{d_1}\left( x \right)} \right)} \right)} \right)} \right)}^{{k \mathord{\left/
 {\vphantom {k {\left( {k + 1} \right)}}} \right.
 \kern-\nulldelimiterspace} {\left( {k + 1} \right)}}}}}}{{{\underline{\xi} _{\varepsilon}}M\left( {{d_1}\left( x \right)} \right)f\left( {\varphi \left( {{\underline{\xi} _{\varepsilon}}M\left( {{d_1}\left( x \right)} \right)} \right)} \right)}}}\right)\\
 \times\sigma_{j-1}\left( \frac{{{\rho _1}\left( {\overline x } \right)}}{{1 - d\left( x \right){\rho _1}\left( {\overline x } \right)}},\cdots,\frac{{{\rho _{n-1}}\left( {\overline x } \right)}}{{1 - d\left( x \right){\rho _{n-1}}\left( {\overline x } \right)}}\right)\\
 +\frac{{M\left( {{d_1}\left( x \right)} \right)}}{{{m}\left( {{d_1}\left( x \right)} \right)}}\frac{{{{\left( {\left( {k + 1} \right)F\left( {\varphi \left( {{\underline{\xi} _{\varepsilon}}M\left( {{d_1}\left( x \right)} \right)} \right)} \right)} \right)}^{{k \mathord{\left/
 {\vphantom {k {\left( {k + 1} \right)}}} \right.
 \kern-\nulldelimiterspace} {\left( {k + 1} \right)}}}}}}{{{\underline{\xi} _{\varepsilon}}M\left( {{d_1}\left( x \right)} \right)f\left( {\varphi \left( {{\underline{\xi} _{\varepsilon}}M\left( {{d_1}\left( x \right)} \right)} \right)} \right)}}\sigma_{j}\left( \frac{{{\rho _1}\left( {\overline x } \right)}}{{1 - d\left( x \right){\rho _1}\left( {\overline x } \right)}},\cdots,\frac{{{\rho _{n-1}}\left( {\overline x } \right)}}{{1 - d\left( x \right){\rho _{n-1}}\left( {\overline x } \right)}}\right)\\
 > 0~~\mbox{in}~~\Omega_{2\delta_{\varepsilon}}/\overline{\Omega}_{\sigma}
\end{array}
\end{equation}
as $\delta_{\varepsilon}>0$ being sufficiently small.

Actually, since, by Lemma 3.1, Lemma 3.2$(i_2)$ and (\ref{CfCm}),
\begin{equation}\label{implim1}
\mathop {\lim }\limits_{\mathop {{x} \in\Omega }\limits_{d(x) \to 0}}\frac{M(d(x))}{m(d(x))}\frac{{{{\left( {\left( {k + 1} \right)F\left( {\varphi \left({\underline{\xi} _{\varepsilon}} {M\left( {d\left( x \right)} \right)} \right)} \right)} \right)}^{{k \mathord{\left/
 {\vphantom {k {(k + 1)}}} \right.
 \kern-\nulldelimiterspace} {(k + 1)}}}}}}{{{\underline{\xi} _{\varepsilon}}M\left( {d\left( x \right)} \right)f\left( {\varphi \left({\underline{\xi} _{\varepsilon}} {M\left( {d\left( x \right)} \right)} \right)} \right)}} = 0
\end{equation}
and
\begin{equation}\label{implim2}
 \begin{array}{l}
 1-\mathop {\lim }\limits_{\mathop {{x} \in\Omega }\limits_{d(x) \to 0}}{ \frac{{M\left( {{d}\left( x \right)} \right)m'\left( {{d}\left( x \right)} \right)}}{{{m^2}\left( {{d}\left( x \right)} \right)}}\frac{{{{\left( {\left( {k + 1} \right)F\left( {\varphi \left( {{\underline{\xi} _{\varepsilon}}M\left( {{d}\left( x \right)} \right)} \right)} \right)} \right)}^{{k \mathord{\left/
 {\vphantom {k {\left( {k + 1} \right)}}} \right.
 \kern-\nulldelimiterspace} {\left( {k + 1} \right)}}}}}}{{{\underline{\xi} _{\varepsilon}}M\left( {{d}\left( x \right)} \right)f\left( {\varphi \left( {{\underline{\xi} _{\varepsilon}}M\left( {{d}\left( x \right)} \right)} \right)} \right)}}}\\
 =1-\frac{1-C_m}{C_f}>0,
\end{array}
\end{equation}
we have, for sufficiently small $\delta_{\varepsilon}>0$,
\begin{equation}\label{imp1}
\frac{M(d(x))}{m(d(x))}\frac{{{{\left( {\left( {k + 1} \right)F\left( {\varphi \left({\underline{\xi} _{\varepsilon}} {M\left( {d\left( x \right)} \right)} \right)} \right)} \right)}^{{k \mathord{\left/
 {\vphantom {k {(k + 1)}}} \right.
 \kern-\nulldelimiterspace} {(k + 1)}}}}}}{{{\underline{\xi} _{\varepsilon}}M\left( {d\left( x \right)} \right)f\left( {\varphi \left({\underline{\xi} _{\varepsilon}} {M\left( {d\left( x \right)} \right)} \right)} \right)}}~~\mbox{is small enough in}~\Omega_{2\delta_{\varepsilon}}
\end{equation}
and
\begin{equation}\label{imp2}
1-{ \frac{{M\left( {{d}\left( x \right)} \right)m'\left( {{d}\left( x \right)} \right)}}{{{m^2}\left( {{d}\left( x \right)} \right)}}\frac{{{{\left( {\left( {k + 1} \right)F\left( {\varphi \left( {{\underline{\xi} _{\varepsilon}}M\left( {{d}\left( x \right)} \right)} \right)} \right)} \right)}^{{k \mathord{\left/
 {\vphantom {k {\left( {k + 1} \right)}}} \right.
 \kern-\nulldelimiterspace} {\left( {k + 1} \right)}}}}}}{{{\underline{\xi} _{\varepsilon}}M\left( {{d}\left( x \right)} \right)f\left( {\varphi \left( {{\underline{\xi} _{\varepsilon}}M\left( {{d}\left( x \right)} \right)} \right)} \right)}}}> 0 \ \ \mbox{in}\ \Omega_{2\delta_{\varepsilon}}.
\end{equation}
Since,
\begin{equation}\label{d1(x)}
0<d_1(x)<2\delta_{\varepsilon}-\sigma~~\mbox{in}~\Omega_{2\delta_{\varepsilon}}/\overline{\Omega}_{\sigma},
\end{equation}
we see that (\ref{imp1}) and (\ref{imp2}) still holds with $d(x)$ being replaced by $d_1(x)$.
Therefore, we obtain (\ref{impud}) from $\mathbf{(a_4)}$ and $\mathbf{(a_6)}$.
By (\ref{superderi}), we have (\ref{superadmi}).

~

Next, we show (\ref{supersolution}).

Since $m$ is nondecreasing, we have, by $\mathbf{(a_3)}$ and (\ref{d1xd2x}),
$$
(\underline{b}-\varepsilon)m^{k+1}(d_1(x))\le(\underline{b}-\varepsilon)m^{k+1}(d(x))\le b(x)~~\mbox{in}~~\Omega_{2\delta_{\varepsilon}}/\overline{\Omega}_{\sigma}.
$$
Hence, to show (\ref{supersolution}), we only need to prove
\begin{equation}\label{similarsuper}
S_k({D^2}{{\overline u} _\varepsilon }\left( x \right)) \leq \left( {{\underline{b}} - \varepsilon } \right){m^{k + 1}}\left( {{d_1}\left( x \right)} \right)f\left( {  {{{\overline u} }_\varepsilon }\left( x \right)}\right) \ \ \mbox{in}\ \Omega_{2\delta_\varepsilon}/\overline{\Omega}_{\sigma}.
\end{equation}

By (\ref{superderi}),
\begin{equation}\label{superderi-simibf}
\begin{array}{l}
S_k({D^2}{{\overline u} _\varepsilon }\left( x \right)) - \left( {{\underline{b}} - \varepsilon } \right){m^{k + 1}}\left( {{d_1}\left( x \right)} \right)f\left( {  {{\overline u }_\varepsilon }\left( x \right)} \right)\\
={\underline{\xi}^{k+1}_{\varepsilon}}m^{k+1}(d_1(x))f(\varphi(\underline{\xi} _{\varepsilon}M(d_1(x))))\\
~~\times\bigg[\left(1-{ \frac{{M\left( {{d_1}\left( x \right)} \right)m'\left( {{d_1}\left( x \right)} \right)}}{{{m^2}\left( {{d_1}\left( x \right)} \right)}}\frac{{{{\left( {\left( {k + 1} \right)F\left( {\varphi \left( {{\underline{\xi} _{\varepsilon}}M\left( {{d_1}\left( x \right)} \right)} \right)} \right)} \right)}^{{k \mathord{\left/
 {\vphantom {k {\left( {k + 1} \right)}}} \right.
 \kern-\nulldelimiterspace} {\left( {k + 1} \right)}}}}}}{{{\underline{\xi} _{\varepsilon}}M\left( {{d_1}\left( x \right)} \right)f\left( {\varphi \left( {{\underline{\xi} _{\varepsilon}}M\left( {{d_1}\left( x \right)} \right)} \right)} \right)}}}\right)\\
 ~~~~~~\times\sigma_{k-1}\left( \frac{{{\rho _1}\left( {\overline x } \right)}}{{1 - d\left( x \right){\rho _1}\left( {\overline x } \right)}},\cdots,\frac{{{\rho _{n-1}}\left( {\overline x } \right)}}{{1 - d\left( x \right){\rho _{n-1}}\left( {\overline x } \right)}}\right)\\
 ~~~~~~+\frac{{M\left( {{d_1}\left( x \right)} \right)}}{{{m}\left( {{d_1}\left( x \right)} \right)}}\frac{{{{\left( {\left( {k + 1} \right)F\left( {\varphi \left( {{\underline{\xi} _{\varepsilon}}M\left( {{d_1}\left( x \right)} \right)} \right)} \right)} \right)}^{{k \mathord{\left/
 {\vphantom {k {\left( {k + 1} \right)}}} \right.
 \kern-\nulldelimiterspace} {\left( {k + 1} \right)}}}}}}{{{\underline{\xi} _{\varepsilon}}M\left( {{d_1}\left( x \right)} \right)f\left( {\varphi \left( {{\underline{\xi} _{\varepsilon}}M\left( {{d_1}\left( x \right)} \right)} \right)} \right)}}\sigma_{k}\left( \frac{{{\rho _1}\left( {\overline x } \right)}}{{1 - d\left( x \right){\rho _1}\left( {\overline x } \right)}},\cdots,\frac{{{\rho _{n-1}}\left( {\overline x } \right)}}{{1 - d\left( x \right){\rho _{n-1}}\left( {\overline x } \right)}}\right)\bigg]\\
 \ \ - ({\underline{b}} - \varepsilon ){m^{k + 1}}\left( {{d_1}\left( x \right)} \right)f\left( {\varphi \left( {{\underline{\xi} _{\varepsilon}}M\left( {{d_1}\left( x \right)} \right)} \right)} \right).
\end{array}
\end{equation}
Then to show (\ref{similarsuper}), we only need to prove
\begin{equation}\label{impineqforsuper}
\begin{array}{l}
{\underline{\xi}^{k+1}_{\varepsilon}}\bigg[\left(1-{ \frac{{M\left( {{d_1}\left( x \right)} \right)m'\left( {{d_1}\left( x \right)} \right)}}{{{m^2}\left( {{d_1}\left( x \right)} \right)}}\frac{{{{\left( {\left( {k + 1} \right)F\left( {\varphi \left( {{\underline{\xi} _{\varepsilon}}M\left( {{d_1}\left( x \right)} \right)} \right)} \right)} \right)}^{{k \mathord{\left/
 {\vphantom {k {\left( {k + 1} \right)}}} \right.
 \kern-\nulldelimiterspace} {\left( {k + 1} \right)}}}}}}{{{\underline{\xi} _{\varepsilon}}M\left( {{d_1}\left( x \right)} \right)f\left( {\varphi \left( {{\underline{\xi} _{\varepsilon}}M\left( {{d_1}\left( x \right)} \right)} \right)} \right)}}}\right)\\
 ~~~~~~\times\sigma_{k-1}\left( \frac{{{\rho _1}\left( {\overline x } \right)}}{{1 - d\left( x \right){\rho _1}\left( {\overline x } \right)}},\cdots,\frac{{{\rho _{n-1}}\left( {\overline x } \right)}}{{1 - d\left( x \right){\rho _{n-1}}\left( {\overline x } \right)}}\right)\\
 ~~~~~~+\frac{{M\left( {{d_1}\left( x \right)} \right)}}{{{m}\left( {{d_1}\left( x \right)} \right)}}\frac{{{{\left( {\left( {k + 1} \right)F\left( {\varphi \left( {{\underline{\xi} _{\varepsilon}}M\left( {{d_1}\left( x \right)} \right)} \right)} \right)} \right)}^{{k \mathord{\left/
 {\vphantom {k {\left( {k + 1} \right)}}} \right.
 \kern-\nulldelimiterspace} {\left( {k + 1} \right)}}}}}}{{{\underline{\xi} _{\varepsilon}}M\left( {{d_1}\left( x \right)} \right)f\left( {\varphi \left( {{\underline{\xi} _{\varepsilon}}M\left( {{d_1}\left( x \right)} \right)} \right)} \right)}}\sigma_{k}\left( \frac{{{\rho _1}\left( {\overline x } \right)}}{{1 - d\left( x \right){\rho _1}\left( {\overline x } \right)}},\cdots,\frac{{{\rho _{n-1}}\left( {\overline x } \right)}}{{1 - d\left( x \right){\rho _{n-1}}\left( {\overline x } \right)}}\right)\bigg]\\
 - ({\underline{b}} - \varepsilon )\le 0~~\mbox{in}~~\Omega_{2\delta_{\varepsilon}}/\overline{\Omega}_{\sigma}
\end{array}
\end{equation}
for sufficiently small $\delta_{\varepsilon}>0$.

Since, by (\ref{LConst}),
\begin{equation}\label{lconst}
\underline{\xi} _{\varepsilon}^{k + 1}\big[(1+\varepsilon){L_0}(1 - C_f^{ - 1}(1 - {C_m}))\big] - ({\underline{b}} - \varepsilon ) =  - \varepsilon,
\end{equation}
we see from (\ref{implim1}), (\ref{implim2}), (\ref{d1(x)}) and $\mathbf{(a_6)}$ that
\begin{equation}\label{iq}
\begin{array}{l}
\underline{\xi}_{\varepsilon}^{k+1}\Bigg[(1+\varepsilon)L_0\left(1-\frac{M(d_1(x))m^\prime(d_1(x))}{{m^2(d_1(x))}}\frac{{{{\left( {\left( {k + 1} \right)F\left( {\varphi \left( \underline{\xi}_{\varepsilon}{M\left( {d_1\left( x \right)} \right)} \right)} \right)} \right)}^{{k \mathord{\left/
 {\vphantom {k {(k + 1)}}} \right.
 \kern-\nulldelimiterspace} {(k + 1)}}}}}}{{\underline{\xi} _{\varepsilon}M\left( {d_1\left( x \right)} \right)f\left( {\varphi \left( \underline{\xi} _{\varepsilon}{M\left( {d_1\left( x \right)} \right)} \right)} \right)}}\right)\\
 ~~~~~~+ \frac{M(d_1(x))}{{m(d_1(x))}}\frac{{{{\left( {\left( {k + 1} \right)F\left( {\varphi \left( \underline{\xi}_{\varepsilon}{M\left( {d_1\left( x \right)} \right)} \right)} \right)} \right)}^{{k \mathord{\left/
 {\vphantom {k {(k + 1)}}} \right.
 \kern-\nulldelimiterspace} {(k + 1)}}}}}}{{\underline{\xi} _{\varepsilon}M\left( {d_1\left( x \right)} \right)f\left( {\varphi \left( \underline{\xi} _{\varepsilon}{M\left( {d_1\left( x \right)} \right)} \right)} \right)}}\sigma_k\left( \frac{{{\rho _1}\left( {\overline x } \right)}}{{1 - d\left( x \right){\rho _1}\left( {\overline x } \right)}},\cdots,\frac{{{\rho _{n-1}}\left( {\overline x } \right)}}{{1 - d\left( x \right){\rho _{n-1}}\left( {\overline x } \right)}}\right)\Bigg]\\
-(\underline{b}-\varepsilon)\le0 ~~\mbox{in}~~\Omega_{2\delta_{\varepsilon}}/\overline{\Omega}_{\sigma}
\end{array}
\end{equation}
for sufficiently small $\delta_{\varepsilon}>0$. This and $\mathbf{(a_5)}$ imply (\ref{impineqforsuper}). Therefore, by (\ref{superderi-simibf}), we obtain (\ref{similarsuper}).

~

\textbf{Step 2.} We prove that ${\underline{u}}_{\varepsilon}$ is $k$-admissible and
\begin{equation}\label{subsolution}
S_k({D^2}{{\underline u} _\varepsilon }\left( x \right)) \geq b(x)f\left( {  {{{\underline u} }_\varepsilon }\left( x \right)}\right)\ \ \mbox{in}\ \Omega_{2\delta_\varepsilon-\sigma}
\end{equation}
as $\delta_{\varepsilon}$ being sufficiently small.

The proof is similar to \textbf{Step 1}.

First, we show that ${\underline{u}}_{\varepsilon}$ is a $k$-admissible function in $\Omega_{2\delta_{\varepsilon}-\sigma}$. That is, for $1\le j\le k$,
\begin{equation}\label{subadmi}
S_j(D^2{\underline{u}}_{\varepsilon})>0~~\mbox{in}~~\Omega_{2\delta_{\varepsilon}-\sigma}.
\end{equation}

As in the proof of \textbf{Step 1}, by straightforward computations, we derive that for $1\le j\le k$,
\begin{equation}\label{subderi}
\begin{array}{l}
S_j(D^2{\underline{u}}_{\varepsilon})\\
={\overline{\xi}^{j+1}_{\varepsilon}}m^{j+1}(d_2(x))f(\varphi(\overline{\xi} _{\varepsilon}M(d_2(x)))){{{\left( {\left( {k + 1} \right)F\left( {\varphi \left( {{\overline{\xi} _{\varepsilon}}M\left( {{d_2}\left( x \right)} \right)} \right)} \right)} \right)}^{{(j-k) \mathord{\left/
 {\vphantom {(j-k) {\left( {k + 1} \right)}}} \right.
 \kern-\nulldelimiterspace} {\left( {k + 1} \right)}}}}}\\
~~\times\bigg[\left(1-{ \frac{{M\left( {{d_2}\left( x \right)} \right)m'\left( {{d_2}\left( x \right)} \right)}}{{{m^2}\left( {{d_2}\left( x \right)} \right)}}\frac{{{{\left( {\left( {k + 1} \right)F\left( {\varphi \left( {{\overline{\xi} _{\varepsilon}}M\left( {{d_2}\left( x \right)} \right)} \right)} \right)} \right)}^{{k \mathord{\left/
 {\vphantom {k {\left( {k + 1} \right)}}} \right.
 \kern-\nulldelimiterspace} {\left( {k + 1} \right)}}}}}}{{{\overline{\xi} _{\varepsilon}}M\left( {{d_2}\left( x \right)} \right)f\left( {\varphi \left( {{\overline{\xi} _{\varepsilon}}M\left( {{d_2}\left( x \right)} \right)} \right)} \right)}}}\right)\\
 ~~~~~~\times\sigma_{j-1}\left( \frac{{{\rho _1}\left( {\overline x } \right)}}{{1 - d\left( x \right){\rho _1}\left( {\overline x } \right)}},\cdots,\frac{{{\rho _{n-1}}\left( {\overline x } \right)}}{{1 - d\left( x \right){\rho _{n-1}}\left( {\overline x } \right)}}\right)\\
 ~~~~~~+\frac{{M\left( {{d_2}\left( x \right)} \right)}}{{{m}\left( {{d_2}\left( x \right)} \right)}}\frac{{{{\left( {\left( {k + 1} \right)F\left( {\varphi \left( {{\overline{\xi} _{\varepsilon}}M\left( {{d_2}\left( x \right)} \right)} \right)} \right)} \right)}^{{k \mathord{\left/
 {\vphantom {k {\left( {k + 1} \right)}}} \right.
 \kern-\nulldelimiterspace} {\left( {k + 1} \right)}}}}}}{{{\overline{\xi} _{\varepsilon}}M\left( {{d_2}\left( x \right)} \right)f\left( {\varphi \left( {{\overline{\xi} _{\varepsilon}}M\left( {{d_2}\left( x \right)} \right)} \right)} \right)}}\sigma_{j}\left( \frac{{{\rho _1}\left( {\overline x } \right)}}{{1 - d\left( x \right){\rho _1}\left( {\overline x } \right)}},\cdots,\frac{{{\rho _{n-1}}\left( {\overline x } \right)}}{{1 - d\left( x \right){\rho _{n-1}}\left( {\overline x } \right)}}\right)\bigg].
\end{array}
\end{equation}
Since
\begin{equation}\label{d2(x)}
\sigma<d_2(x)<2\delta_{\varepsilon}~~\mbox{in}~\Omega_{2\delta_{\varepsilon}-\sigma},
\end{equation}
by the same argument as in \textbf{Step 1}, we also have, for
$1\le j\le k$,
\begin{equation}\label{impld}
\begin{array}{l}
\left(1-{ \frac{{M\left( {{d_2}\left( x \right)} \right)m'\left( {{d_2}\left( x \right)} \right)}}{{{m^2}\left( {{d_2}\left( x \right)} \right)}}\frac{{{{\left( {\left( {k + 1} \right)F\left( {\varphi \left( {{\overline{\xi} _{\varepsilon}}M\left( {{d_2}\left( x \right)} \right)} \right)} \right)} \right)}^{{k \mathord{\left/
 {\vphantom {k {\left( {k + 1} \right)}}} \right.
 \kern-\nulldelimiterspace} {\left( {k + 1} \right)}}}}}}{{{\overline{\xi} _{\varepsilon}}M\left( {{d_2}\left( x \right)} \right)f\left( {\varphi \left( {{\overline{\xi} _{\varepsilon}}M\left( {{d_2}\left( x \right)} \right)} \right)} \right)}}}\right)\\
 \times\sigma_{j-1}\left( \frac{{{\rho _1}\left( {\overline x } \right)}}{{1 - d\left( x \right){\rho _1}\left( {\overline x } \right)}},\cdots,\frac{{{\rho _{n-1}}\left( {\overline x } \right)}}{{1 - d\left( x \right){\rho _{n-1}}\left( {\overline x } \right)}}\right)\\
 +\frac{{M\left( {{d_2}\left( x \right)} \right)}}{{{m}\left( {{d_2}\left( x \right)} \right)}}\frac{{{{\left( {\left( {k + 1} \right)F\left( {\varphi \left( {{\overline{\xi} _{\varepsilon}}M\left( {{d_2}\left( x \right)} \right)} \right)} \right)} \right)}^{{k \mathord{\left/
 {\vphantom {k {\left( {k + 1} \right)}}} \right.
 \kern-\nulldelimiterspace} {\left( {k + 1} \right)}}}}}}{{{\overline{\xi} _{\varepsilon}}M\left( {{d_2}\left( x \right)} \right)f\left( {\varphi \left( {{\overline{\xi} _{\varepsilon}}M\left( {{d_2}\left( x \right)} \right)} \right)} \right)}}\sigma_{j}\left( \frac{{{\rho _1}\left( {\overline x } \right)}}{{1 - d\left( x \right){\rho _1}\left( {\overline x } \right)}},\cdots,\frac{{{\rho _{n-1}}\left( {\overline x } \right)}}{{1 - d\left( x \right){\rho _{n-1}}\left( {\overline x } \right)}}\right)\\
 > 0~~\mbox{in}~~\Omega_{2\delta_{\varepsilon}-\sigma}
\end{array}
\end{equation}
as $\delta_{\varepsilon}>0$ being sufficiently small.
This combining with (\ref{subderi}) imply (\ref{subadmi}).

~

Next, we show (\ref{subsolution}).

Since $m$ is nondecreasing, we have, by $\mathbf{(a_3)}$ and (\ref{d1xd2x}),
$$
(\overline{b}+\varepsilon)m^{k+1}(d_2(x))\ge(\overline{b}+\varepsilon)m^{k+1}(d(x))\ge b(x)~~\mbox{in}~~\Omega_{2\delta_{\varepsilon}-\sigma}.
$$
Hence, (\ref{subsolution}) is an easy consequence of
\begin{equation}\label{similarsub}
S_k({D^2}{{\underline u} _\varepsilon }\left( x \right)) \geq \left( {{\overline{b}} + \varepsilon } \right){m^{k + 1}}\left( {{d_2}\left( x \right)} \right)f\left( {  {{{\underline u} }_\varepsilon }\left( x \right)}\right) \ \ \mbox{in}\ \Omega_{2\delta_\varepsilon-\sigma}.
\end{equation}

It follows from (\ref{subderi}) that to show (\ref{similarsub}), we only need to prove
\begin{equation}\label{impineqforsub}
\begin{array}{l}
{\overline{\xi}^{k+1}_{\varepsilon}}\bigg[\left(1-{ \frac{{M\left( {{d_2}\left( x \right)} \right)m'\left( {{d_2}\left( x \right)} \right)}}{{{m^2}\left( {{d_2}\left( x \right)} \right)}}\frac{{{{\left( {\left( {k + 1} \right)F\left( {\varphi \left( {{\overline{\xi} _{\varepsilon}}M\left( {{d_2}\left( x \right)} \right)} \right)} \right)} \right)}^{{k \mathord{\left/
 {\vphantom {k {\left( {k + 1} \right)}}} \right.
 \kern-\nulldelimiterspace} {\left( {k + 1} \right)}}}}}}{{{\overline{\xi} _{\varepsilon}}M\left( {{d_2}\left( x \right)} \right)f\left( {\varphi \left( {{\overline{\xi} _{\varepsilon}}M\left( {{d_2}\left( x \right)} \right)} \right)} \right)}}}\right)\\
 ~~~~~~\times\sigma_{k-1}\left( \frac{{{\rho _1}\left( {\overline x } \right)}}{{1 - d\left( x \right){\rho _1}\left( {\overline x } \right)}},\cdots,\frac{{{\rho _{n-1}}\left( {\overline x } \right)}}{{1 - d\left( x \right){\rho _{n-1}}\left( {\overline x } \right)}}\right)\\
 ~~~~~~+\frac{{M\left( {{d_2}\left( x \right)} \right)}}{{{m}\left( {{d_2}\left( x \right)} \right)}}\frac{{{{\left( {\left( {k + 1} \right)F\left( {\varphi \left( {{\overline{\xi} _{\varepsilon}}M\left( {{d_2}\left( x \right)} \right)} \right)} \right)} \right)}^{{k \mathord{\left/
 {\vphantom {k {\left( {k + 1} \right)}}} \right.
 \kern-\nulldelimiterspace} {\left( {k + 1} \right)}}}}}}{{{\overline{\xi} _{\varepsilon}}M\left( {{d_2}\left( x \right)} \right)f\left( {\varphi \left( {{\overline{\xi} _{\varepsilon}}M\left( {{d_2}\left( x \right)} \right)} \right)} \right)}}\sigma_{k}\left( \frac{{{\rho _1}\left( {\overline x } \right)}}{{1 - d\left( x \right){\rho _1}\left( {\overline x } \right)}},\cdots,\frac{{{\rho _{n-1}}\left( {\overline x } \right)}}{{1 - d\left( x \right){\rho _{n-1}}\left( {\overline x } \right)}}\right)\bigg]\\
 - ({\overline{b}} + \varepsilon )\ge 0~~\mbox{in}~~\Omega_{2\delta_{\varepsilon}-\sigma}.
\end{array}
\end{equation}

Since, by (\ref{UConst}),
\begin{equation}\label{uconst}
\overline{\xi} _{\varepsilon}^{k + 1}\big[(1-\varepsilon){l_0}(1 - C_f^{ - 1}(1 - {C_m}))\big] - ({\overline{b}} + \varepsilon ) =  \varepsilon,
\end{equation}
we see from $(\ref{implim1})$, $(\ref{implim2})$ (note here that $\underline{\xi} _{\varepsilon}$ is replaced by $\overline{\xi} _{\varepsilon}$), (\ref{d2(x)}) and $\mathbf{(a_6)}$ that
\begin{equation}
\begin{array}{l}
\overline{\xi}_{\varepsilon}^{k+1}\Bigg[(1-\varepsilon)l_0\left(1-\frac{M(d_2(x))m^\prime(d_2(x))}{{m^2(d_2(x))}}\frac{{{{\left( {\left( {k + 1} \right)F\left( {\varphi \left( \overline{\xi}_{\varepsilon}{M\left( {d_2\left( x \right)} \right)} \right)} \right)} \right)}^{{k \mathord{\left/
 {\vphantom {k {(k + 1)}}} \right.
 \kern-\nulldelimiterspace} {(k + 1)}}}}}}{{\overline{\xi} _{\varepsilon}M\left( {d_2\left( x \right)} \right)f\left( {\varphi \left( \overline{\xi} _{\varepsilon}{M\left( {d_2\left( x \right)} \right)} \right)} \right)}}\right)\\
 ~~~~~~+ \frac{M(d_2(x))}{{m(d_2(x))}}\frac{{{{\left( {\left( {k + 1} \right)F\left( {\varphi \left( \overline{\xi}_{\varepsilon}{M\left( {d_2\left( x \right)} \right)} \right)} \right)} \right)}^{{k \mathord{\left/
 {\vphantom {k {(k + 1)}}} \right.
 \kern-\nulldelimiterspace} {(k + 1)}}}}}}{{\overline{\xi} _{\varepsilon}M\left( {d_2\left( x \right)} \right)f\left( {\varphi \left( \overline{\xi} _{\varepsilon}{M\left( {d_2\left( x \right)} \right)} \right)} \right)}}\sigma_k\left( \frac{{{\rho _1}\left( {\overline x } \right)}}{{1 - d\left( x \right){\rho _1}\left( {\overline x } \right)}},\cdots,\frac{{{\rho _{n-1}}\left( {\overline x } \right)}}{{1 - d\left( x \right){\rho _{n-1}}\left( {\overline x } \right)}}\right)\Bigg]\\
-(\overline{b}+\varepsilon)\ge 0 ~~\mbox{in}~~\Omega_{2\delta_{\varepsilon}}
\end{array}
\end{equation}
for sufficiently small $\delta_{\varepsilon}>0$. This combining with $\mathbf{(a_5)}$ imply (\ref{impineqforsub}). Consequently, we have (\ref{similarsub}) and then (\ref{subsolution}).

~

\textbf{Step 3.} We show (\ref{MBB}).

Let $u\in C(\Omega)$ be a viscosity solution of (1.1) and let $T>0$ (depending on $\delta_{\varepsilon}$) sufficiently large such that for any $0<\sigma<\delta_{\varepsilon}$,
\begin{equation}\label{cb1}
u\leq {\overline{u}}_{\varepsilon}+T\ \ \mbox{on}\ \Lambda_1=\{x\in\Omega: d(x)=2\delta_{\varepsilon}\}
\end{equation}
and
\begin{equation}\label{acb1}
{\underline{u}}_{\varepsilon}\leq u+T\ \ \mbox{on}\ \Lambda_2=\{x\in\Omega: d(x)=2\delta_{\varepsilon}-\sigma\}.
\end{equation}
We observe that
\begin{equation}\label{cb2}
u\leq {\overline{u}}_{\varepsilon}+T=\infty \ \ \mbox{on}\ \Lambda_2=\{x\in\Omega: d(x)=\sigma\}.
\end{equation}
and
\begin{equation}\label{acb2}
{\underline{u}}_{\varepsilon}\le u+T=\infty \ \ \mbox{on}\ \partial\Omega.
\end{equation}

Since, $f$ is nondecreasing and ${\overline{u}}_{\varepsilon}$ satisfies (\ref{supersolution}), we have
$$
S_k(D^2({\overline{u}}_{\varepsilon}+T))= S_k( D^2{\overline{u}}_{\varepsilon})\leq b(x)f({\overline{u}}_{\varepsilon})\leq b(x)f({\overline{u}}_{\varepsilon}+T)\ \ \mbox{in}\ \Omega_{2\delta_{\varepsilon}}/\Omega_{\sigma}.
$$
Note that in the viscosity sense,
$$
S_k( D^2u))= b(x)f(u)\ \ \mbox{in}\ \Omega.
$$
Therefore, by Lemma 2.1, we deduce from (\ref{cb1}) and  (\ref{cb2}) that
\begin{equation}\label{cp}
u\leq {\overline{u}}_{\varepsilon}+T\ \ \mbox{in}\ \Omega_{2\delta_{\varepsilon}}/\Omega_{\sigma}.
\end{equation}
Similarly, since in the viscosity sense
$$
S_k(D^2(u+T))= S_k( D^2u)= b(x)f(u)\leq b(x)f(u+T)\ \ \mbox{in}\ \Omega_{2\delta_{\varepsilon}-\sigma},
$$
and $\underline{u}_{\varepsilon}$ satisfies (\ref{subsolution}), we deduce from (\ref{acb1}) and (\ref{acb2}) that
\begin{equation}\label{acp}
{\underline{u}}_{\varepsilon}\leq {u+T}\ \ \mbox{in}\ \Omega_{2\delta_{\varepsilon}-\sigma}.
\end{equation}

Substituting (\ref{supersub}) into (\ref{cp}) and (\ref{acp}) respectively, we have
$$
\frac{u}{\varphi(\underline{\xi} _{\varepsilon}M(d_1(x)))}\le 1+\frac{T}{\varphi(\underline{\xi} _{\varepsilon}M(d_1(x)))}\ \ \mbox{in}\ \Omega_{2\delta_{\varepsilon}}/\Omega_{\sigma}
$$
and
$$
 1-\frac{T}{\varphi(\overline{\xi} _{\varepsilon}M(d_2(x)))} \le \frac{u}{\varphi(\overline{\xi} _{\varepsilon}M(d_2(x)))} \ \ \mbox{in}\ \Omega_{2\delta_{\varepsilon}-\sigma}.
$$
Let $\sigma\rightarrow 0$,
$$
\frac{u}{\varphi(\underline{\xi} _{\varepsilon}M(d(x)))}\le 1+\frac{T}{\varphi(\underline{\xi} _{\varepsilon}M(d(x)))}\ \ \mbox{in}\ \Omega_{2\delta_{\varepsilon}}
$$
and
$$
 1-\frac{T}{\varphi(\overline{\xi} _{\varepsilon}M(d(x)))}\le \frac{u}{\varphi(\overline{\xi} _{\varepsilon}M(d(x)))} \ \ \mbox{in}\ \Omega_{2\delta_{\varepsilon}}.
$$
Note that
$$
\varphi(\underline{\xi} _{\varepsilon}M(d(x)))=\infty\ \ \mbox{on}\ \partial\Omega
$$
and
$$
\varphi(\overline{\xi} _{\varepsilon}M(d(x)))=\infty\ \ \mbox{on}\ \partial\Omega.
$$
We obtain
$$ \mathop {\lim \sup}\limits_{\mathop {{x} \in\Omega }\limits_{d(x) \to 0}}  \frac{{u(x)}}{{\varphi ({\underline{\xi} _{\varepsilon}}M(d(x)))}}\le 1
$$
and
$$
1\le\mathop {\liminf}\limits_{\mathop {{x} \in\Omega }\limits_{d(x) \to 0}}  \frac{{u(x)}}{{\varphi ({\overline{\xi}_{\varepsilon} }M(d(x)))}}.
$$
Let $\varepsilon\rightarrow 0$ and then we conclude (1.6).

The proof of Theorem 1.1 is complete.

\qed

\section*{References}

\bibliographystyle{plain}

\end{document}